\documentclass{amsart}

\newtheorem{theorem}{Theorem}
\newtheorem{proposition}{Proposition}

\theoremstyle{definition}

\theoremstyle{remark}
\newtheorem{remark}{Remark}

\numberwithin{equation}{section}

\newcommand{\ba}{\begin{array}}
\newcommand{\be}{\begin{equation}}
\newcommand{\ea}{\end{array}}
\newcommand{\ee}[1]{\label{#1}\end{equation}}

\newcommand{\Cal}[1]{{\mathcal {#1}}}

\newcommand{\co}{{\mathfrak{o}}}

\newcommand{\cs}{{\mathfrak{s}}}

\begin{document}

\author{Sergio Console}
\address{Dipartimento di Matematica
Universit\`a di Torino,
via Carlo Alberto 10,
10123 Torino, Italy}
\email{sergio.console@unito.it}
\thanks{The first author was partially
supported by GNSAGA of INdAM, MIUR of Italy, CONICET, Secyt-UNC and CIEM of Argentina}

\author{Carlos Olmos}
\address{FaMAF, Universidad Nacional de C\'ordoba,
Ciudad Universitaria,
5000 C\'ordoba, Argentina}
\email{olmos@mate.uncor.edu}
\thanks{The second author was supported by Universidad Nacional de C\'ordoba and CONICET,
partially supported by   Antorchas, ANCyT, Secyt-UNC and CIEM}




\title[Curvature invariants, Killing fields, connections and cohomogeneity]{Curvature invariants, Killing vector fields, connections and cohomogeneity}

\begin{abstract}
A direct, bundle-theoretic method for defining and extending local isometries out of curvature data is developed. As a by-product, conceptual direct proofs of a classical result of Singer and a recent result of the authors are derived. 
\end{abstract}

\maketitle

A classical result of I. M. Singer \cite{Si} states that a Riemannian manifold is locally homogeneous if and only if its Riemannian curvature tensor together with its covariant derivatives up to some index $k+1$ are independent of the point (the integer $k$ is called the Singer invariant). More precisely

\begin{theorem}[Singer \cite{Si}] \label{singer}
Let $M$ be a Riemannian manifold. Then $M$ is locally homogeneous if and only if for any $p, q \in M$ there is a linear isometry $F: T_p M \to T_qM$ such that 
$$
F^* \nabla^s R_q= \nabla^s R_p\, ,
$$
for any $s \leq k+1$. 
\end{theorem}

An alternate proof with a more direct approach was given in \cite{NT}.

Out of the curvature tensor and its covariant derivatives one can construct scalar invariants, like for instance the scalar curvature. In general, any polynomial function in the components of  the curvature tensor and its covariant derivatives which does not depend on the choice of the orthonormal basis at the tangent space of each point is a \textit{a scalar Weyl invariant} or  \textit{a scalar curvature invariant}.
By Weyl  theory of invariants, 
a scalar Weyl  invariant is a linear combination of complete traces of tensors
$\langle \nabla^{m_1} R, \, \rangle \dots \langle \nabla^{m_\ell} R, \, \rangle$, ($m_1, \dots m_s \geq 0$, $\nabla^0 R=R$).

Pr\"ufer, Tricerri and Vanhecke studied the interplay among local homogeneity and these curvature invariants. Using Singer's Theorem, 
they got the following 

\begin{theorem}[Pr\"ufer, Tricerri and Vanhecke \cite{PTV}] \label{ptv}
Let $M$ be an $n$-dimensional Riemannian manifold. Then $M$ is locally homogeneous if and only if all scalar Weyl invariants of order $s$ with $s \leq \frac{n(n-1)}{2}$ are constant. 
\end{theorem}

More in general for a non-homogeneous Riemannian manifold one can look at the regular level sets of  scalar Weyl invariants. 
In a recent paper \cite{CO}, we proved the following theorem, which generalizes the above results

\begin{theorem}[\cite{CO}] \label {main-co} The cohomogeneity of a Riemannian manifold $M$ (with respect to the full isometry group) coincides locally with the codimension of the foliation by
regular level sets of the scalar Weyl invariants.
\end {theorem}

The key point in the proof was extending Killing vector fields in any level set to $M$ (at least locally). In the present paper we use a direct method, identifying Killing fields with parallel sections of a bundle on $M$ (namely the one whose fiber at $p$ is $E_p = T_pM\oplus \cs\co (T_pM)$, where $\cs\co (T_pM)$ are the skew-symmetric endomorphisms of $T_pM$) endowed with a connection $\tilde \nabla$ (induced by the Levi-Civita connection of $M$). This local characterization of Killing vector fields
as parallel sections of a vector bundle with connection goes back to Kostant \cite{Ko} and it is used in theoretical physics (see for instance \cite{BFP}). In the setting of Theorem~\ref{main-co},  it turns out that Killing fields in level sets correspond to a subbundle $\bar E$ of $E$, which we prove to be parallel and flat. This yields a conceptual and direct proof of Theorem~\ref{main-co} and Singer's Theorem, 
which does not use homogeneous structures.

 \section{The connection}

Let $M$ be a Riemannian manifold and let $E \overset {\pi}{\rightarrow} M$ be the metric bundle over $M$ whose
fibers are $E_p = T_pM\oplus \cs\co (T_pM)$.
The bundle $E$ is
canonically isomorphic with $TM\oplus \Lambda ^2(M)$. We endow $E$ with the connection $\tilde \nabla$, induced
by the Levi-Civita connection $\nabla$ on $M$. Namely, if $(v,B)$ is a section of $E$ (i.e., $v$ is a vector
field and $B$ is skew-symmetric  tensor field of type $(1,1)$ on $M$)  then
$$\tilde \nabla _X (v,B) = (\nabla _X v - B.X,  \nabla _X B - R_{X,v})$$
where $R$ is the curvature tensor on $M$ and, as usual, $(\nabla _X B).Y = \nabla _X (B.Y) - B. \nabla _X Y$.

The {\it canonical lift} of vector field $Z$ on $M$ is the section $\hat Z$ of $E$ given by
$$ \hat Z (p) =  (Z(p), [(\nabla Z)_p]^{\text {skew}})$$
where the upper script  ``skew'' denotes the skew symmetric part of the endomorphisms  $(\nabla Z)_p$ of $T_pM$.

The following result is well known and elementary to show. For the sake of self completeness we include the
proof (cf. also \cite[Section 3.5.2]{BFP}).  \medskip

\begin{proposition}The canonical lift gives an isomorphism between the set $\Cal K (M)$
of Killing fields on $M$ and the parallel sections of $E$ with respect to $\tilde \nabla$
\end{proposition}

\begin{proof}  A vector field  $Z$ on $M$ is a Killing field if and only if $(\nabla Z)_p$ is skew symmetric.
Observe that in this case $Z$ satisfies the affine Jacobi equation.
$$ \nabla _X (\nabla Z) - R_{X,Z}=0$$
for all $X$.  This equation is derived from the fact that the associated flow to $Z$ preserves the Levi-Civita
connection (by making use of the fact that $\nabla$ is torsion free and the first Bianchi identity).

So, if $Z$ is Killing then $\hat Z$ is a parallel section of $E$, since $\nabla Z$ is skew-symmetric.
Conversely, if $(v,B)$ is a parallel section of $E$, then the first component of $\tilde \nabla (v, B) = 0$
implies that $\nabla v$ is skew-symmetric and hence $v$ is a Killing field on $M$.
\end{proof}

It is straightforward to compute the curvature tensor of $E$, by making use of the first and the second Bianchi identity (for the first and the second component respectively)

\begin{proposition} The curvature tensor of $E$, with respect to the connection $\tilde \nabla$ is given
by
$$\tilde R_{X,Y}(v, B) = (0, (\nabla _v R)_{X,Y} - (B.R)_{X,Y})$$
where $B$ acts on $R$ as a derivation.
\end{proposition}

\medskip

\begin{proposition} Let $T$ be a given tensor on $M$ and let $(v,B)$ be a section of $E$ that satisfies
the equation
$$\nabla _vT = B.T$$ We have that  $\tilde \nabla _X (v,B)$  satisfies this
equation, for all vector fields $X$ on $M$, if and only if $(v,B)$ also satisfies the following equation
$$ \nabla  _v (\nabla T) = B.(\nabla T)$$
\end{proposition}

\begin{proof} It is straightforward and makes use of the so called Ricci identity $\nabla _{X,Y}T- \nabla _{Y,X}T
= R_{X,Y}.T$, where $R_{X,Y}$ acts as a derivation.
\end{proof}
 
  \section{Extension of Killing vector fields}

Let $M$ be a Riemannian manifold. By making $M$ possibly smaller we may assume that $M$ is foliated by the
regular level sets of the scalar Weyl invariants (cf.  \cite[Section 3]{CO}). Given $p,\, q$ in the same level set, let us say $F$, then
there exists a linear isometry  $h: T_pM \to T_qM$ which maps any covariant derivative $(\nabla ^ k R)_q$ to the
same object at $p$, for any $k\geq 0$.

Let $c(t)$ be a  curve in $F$ with $c(0) = p$ and let $\tau _t : T_pM \to T_{c(t)}M$ be the parallel transport
along $c(t)$. Since the parallel transport is a linear isometry (from the corresponding tangent spaces) and by
the previous paragraph one has that $\tau _t ^{-1} (\nabla ^ k R)_{c(t)}$ lies in the same $O(T_pM)$ orbit of
$(\nabla ^ k R)_p$, for all $t$. Derivating this condition at $t=0$ yields that there exists $B \in\cs \co (T_pM)$ such that
$$ (\nabla _v (\nabla ^ k R))_p = B. (\nabla ^ k R)_p$$
where $v = c'(0) \in T_pF$. Observe that $v$ is arbitrary in $T_pF$, since $c(t)$ is an arbitrary curve in $F$.

Let, for $q \in M$, $E^k_q $ be the subspace of $E_q$ which consists of all the pairs $(v,B)$ such that

$$ (\nabla _v (\nabla ^ i R))_q - B. (\nabla ^ i R)_q = 0$$
for all $0\leq i\leq k$. Notice that the projection to the first component maps $E^k_q$ onto $T_qF(q)$, where
$F(q)$ is the level set of the scalar Weyl invariants by $q$. So, $\text {dim\, }E^k_ q \geq r$, where $r$ is
the dimension of $F(q)$. It is clear that there exists $j(q)\leq \text {dim\, } E_q = n + \frac {1}{2}n(n-1)$
such that $\text {dim\, } E^{j(q)} = \text {dim\, } E^{j(q)+ 1}_q$

By making possibly $M$ smaller we may assume that $\text {dim\, }E^{j(q)}_ q$ does not depend on $q$. This
gives rise to a subbundle $\bar E$ of $E$ whose fibers are $E^{j(q)}_ q$. By Proposition 3 we have that $\bar E$
is a parallel subbundle of $E$ and by Proposition 2 one has that it is flat.  Therefore any $(v,B) \in \bar
E_p$, $p$ fixed in $M$, gives rise to a parallel section $(\tilde v, \tilde B)$ of  $\bar E$ and so to a
parallel section of $E$ (we may assume that $M$ is simply connected). By Proposition 1, this parallel section
corresponds to a Killing field on $M$, whose value at $p$ is $v$, which is arbitrary in $T_pF(p)$. This Killing
field must be always tangent to any level set, because the scalar Weyl invariants are preserved by isometries.

This finishes a shorter, conceptual proof of Theorem~\ref{main-co} (and Singer's Theorem).
\bigskip

\begin{remark}[On the pseudo Riemannian case]
Singer's Theorem generalizes to the pseudo Riemannian case (Podest\`a and Spirio, \cite{PS}). Our proof can be applied to this setting and more generally  to any affine connection without torsion. 
\newline
For pseudo Riemannian manifolds, the behavior of the curvature tensor and its covariant derivatives differs from the one of scalar curvature invariants.
\newline
Indeed, there are examples of (not locally  homogeneous) lorentzian manifolds whose scalar curvature invariants vanish (see e.g. \cite{BP, CHP, KM} and the references therein). Therefore, in this case, having constant scalar curvature invariants does not imply local homogeneity. 
\end{remark}

 \bibliographystyle{amsplain}

\end{document}